\documentclass[12pt]{article}
\def\eqnum#1{\eqno (#1)}
\def\fnote#1{\footnote}
\def\blacksquare{\hfill\hbox{\vrule width 4pt height 4pt depth 0pt}}

\begin{document}
\title{The lattice of varieties of semigroups\\
with completely regular square. II\thanks{\textbf{Added in September 2026:} This paper was written in 1991 as a continuation of \cite{11}, but it was never submitted for publication. In the present version, a few typos have been corrected.}}
\author{T.A.Ershova and M.V.Volkov}

\date{}

\maketitle

\begin{abstract}
We solve (modulo groups) the word problem for free semigroup satisfying
$xy = (xy)^n$. This enables us to show that the variety ${\cal {SCR}}_n$ given
by this identity is not equal to the join of the varieties ${\cal{CR}}_n$ and
${\cal{SI}}$ defined by the identities $x = x^n$ and $xy = (xy)^2$
respectively. Therefore the result of [11] that the lattice $L({\cal {SCR}}_n)$
is modular does not follow from the earlier results concerning $L({\cal{CR}}_n)$
and $L({\cal{SI}})$. However the variety ${\cal{SO}}_n$ consisting of all
semigroups from ${\cal {SCR}}_n$ whose idempotents form a subsemigroup turns out
to be equal to the join of the varieties ${\cal O}_n$ of all orthogroups
and ${\cal{SI}}$. We deduce from this a description of the lattice
$L({\cal{SO}}_n)$ modulo groups.
\end{abstract}
\par
\vskip .3cm
\centerline{\bf Introduction and summary}
\vskip .3cm
A semigroup $S$ is said to be a {\it semigroup with completely regular
square\/} if $S^2 =\{st : s,t\in S\}$ happens to be a completely regular semigroup (i.e. a union
of groups). It is easy to see that a semigroup variety consists of semigroups
with completely regular square if and only if, for some $n > 1$, it satisfies
the identity
$$(xy)^n = xy.
\eqnum{1}$$
Let us fix $n$ and denote by ${\cal {SCR}}_n$ the variety defined by (1). In [11],
where we began the study of this variety, it was shown that the lattice
$L({\cal {SCR}}_n)$ of all subvarieties of ${\cal {SCR}}_n$ is modular (in fact, even
arguesian).\footnote{\textbf{Added in September 2026:} V.~Yu.~Shaprynski\v{\i} found a lacuna in the proof of this result in~\cite{11} and showed how it can be filled. This omission does not affect the specific results from~\cite{11} used in the present paper.} Before that, Pastijn [6] has proved that the lattice of all completely regular semigroup varieties is arguesian and Gerhard [3] has
described the lattice of all varieties of semigroups having idempotent square---i.e., satisfying the identity $(xy)^2 = xy$---and found that it is even
distributive [3]. Comparing these results, we immediately encounter the
following natural and important question: does the variety ${\cal {SCR}}_n$
coincide with the join of the completely regular semigroup variety
${\cal{CR}}_n$ determined by the identity $x^n = x$ and the variety ${\cal{SI}}$
of all semigroups having idempotent square? If the answer would be
affirmative, the lattice $L({\cal {SCR}}_n)$ would be equal to a subdirect product
of the lattices $L({\cal{CR}}_n)$ and $L({\cal{SI}})$ and therefore the result of
[11] would follow easily from the results of [6] and [3] mentioned above.
Here we show that the answer is (fortunately) negative. A solution to the word
problem for free ${\cal {SCR}}_n$-semigroups plays an essential role in our proof
of this fact.
\par
On the other hand, considering the subvariety ${\cal{SO}}_n$ of ${\cal {SCR}}_n$
consisting of semigroups with orthogroup square\fnote{1}
{\noindent Let us recall that an {\it orthogroup\/} is a completely
regular semigroup whose idempotents form a subsemigroup.},
we prove that it does coincide with the join of ${\cal{SI}}$ and ${\cal O}_n$
(where ${\cal O}_n$ is the variety of all orthogroups satisfying $x^n = x$).
Since the lattice $L({\cal O}_n)$ is known modulo group varieties (see [9] or [8]), we
are able to deduce from the latter fact a description of the lattice
$L({\cal{SO}}_n)$.
\par
The paper consists of four sections. Section 1 contains some preliminaries.
Being the core  of the paper, Section 2 is
devoted to the solution to the word problem for free ${\cal {SCR}}_n$-semigroups.
In Section 3 we calculate the join of ${\cal{CR}}_n$ and ${\cal{SI}}$ and verify
that it does not coincide with the variety ${\cal {SCR}}_n$. In Section 4 we
consider varieties of semigroups with orthogroup square.
\par
The authors are indebted to  F.Pastijn for his attention to the
paper and many helpful discussions.
\par
\vskip 0.3cm
\centerline{\bf 1. Preliminaries}
\vskip 0.3cm
We use the standard notation and terminology of semigroup theory (see e.g. [1]).
In addition, for $X$ a fixed countably infinite set, the free
semigroup over $X$ is denoted by $X^+$. The identity relation on $X^+$
is denoted by $\equiv $. By $X^*$ we denote $X^+$ with the empty word
adjoined. Let $w\in X^*$. We define the {\it content} of $w$ (denoted
by $c(w)$) as the set of letters of $X$ occurring in the word $w$. By $| w| $ we denote
length of a word $w$. By the
{\it left (right) indicator\/} of $w$ we mean the shortest initial (final)
segment of $w$ whose content coincides with $c(w)$. The left (right)
indicator of $w$ has the form $ux$ where $u\in X^*, x\in X, c(u) =
c(w)\backslash \{x\}$ (resp., $yv$ where $v\in X^*, y\in X, c(v) =
c(w)\backslash \{y\})$; we denote the word $u$ by $s(w)$, the letter $x$
by $\sigma (w)$ (the word $v$ by $d(w)$ and the letter $y$ by $\delta(w))$.
We put also $s^2(w) \equiv s(s(w))$ and more generally $s^{k}(w) \equiv
s(s^{k-1}(w))$, $k\le| c(w)| $.
\par
The fully invariant congruence on $X^+$ corresponding to a variety
${\cal A}$ is denoted by $\rho _{{\cal A}}$. The semigroup
$X^+/\rho _{\cal A}$ represents the free object on $X$ of the variety
${\cal A}$ and is denoted by $F{\cal A}$. The class of a congruence
$\rho $ containing an element $u$ is denoted by $u^{\rho }$. By $var \Sigma $
we denote the variety given by the set $\Sigma $ of identities.
\par
Now we recall some results from [11] which are important in what follows.
\par
\noindent{\bf RESULT 1} [11, Proposition 1.6]. {\it Let} $u,v\in X^+$. {\it Then}
$u\rho _{{\cal {SCR}}_n}v$ {\it implies}
$s(u)\rho _{{\cal {SCR}}_n}s(v), \sigma (u) \equiv \sigma (v), c(u) = c(v)$.
\blacksquare
\par
This fact and its dual allow us  to speak about the indicators and the content of
an element of the semigroup $F = F{\cal {SCR}}_n$. In the following two results
we use these notions
to describe Green's relations on $F$.
\par
\noindent{\bf RESULT 2} [11, Proposition 1.5]. {\it Let} $u,v\in F$. {\it Then}
$u{\cal D}v$ {\it if and only if either} $u = v$ {\it in F or} $u,v\in F^2$
{\it and} $c(u) = c(v)$. \blacksquare
\par
\noindent{\bf RESULT 3} [11, Proposition 1.8]. {\it Let} $u,v\in F$. {\it Then}
$u{\cal R}v$ {\it if and only if either} $u = v$ {\it in F or} $u,v\in F^2$
{\it and they have the same left indicator in F.} \blacksquare
\par
The last result from [11] we need below is an "${\cal SCR}_n$-analogue" of the so-called
Brown lemma.
\par
\noindent{\bf RESULT 4} [11, Corollary 1.9]. {\it Let} $u\in F^2, v\in F$ {\it and}
$c(v) \subseteq c(u).$ {\it Then} $u = (uv)^{n-1}u = u(vu)^{n-1}$.
\blacksquare
\vskip 0.3cm
\centerline{\bf 2. The word problem}
\vskip 0.3cm
Let ${\cal {V}}_n$ be a variety of groups satisfying the identity
$$x^{n-1} = 1
\eqnum{2}$$
for $n\ge 2$.
\par
We denote by ${\cal{SCRV}}_n$
the variety of all semigroups with completely regular square whose maximal
subgroups belong to ${\cal {V}}_n$. Our aim is to solve the word problem for its
free semigroup over $X$ modulo a solution to the word problem for $F{\cal {V}}_n$.
For brevity, throughout this section, we denote $F{\cal{SCRV}}_n$ by $T$, $\rho
_{{\cal{SCRV}}_n}$
by
$\tau $ and $\rho _{{\cal {SCR}}_n}$ by $\rho $ .
\par
An effective solution to the word problem for free objects in the variety of
all completely regular semigroups whose subgroups are in ${\cal {V}}_n$ was found
in [7] (and in [5] for the case where ${\cal {V}}_n$ is the variety of all groups
satisfying (2)). We follow the approach developed in these papers.
Let us mention also that the word problem for the variety ${\cal{SI}}$ was
solved by Gerhard in [2].
\par
Let us start with constructing an identity basis for the variety ${\cal{SCRV}}_n$.
If $u\in X^+$ and $c(u) \subseteq \{x_{1},\ldots ,x_{k}\}$, we
denote the result of substituting $g_{i}\in X^+$ for $x_{i}
(i=1,\ldots ,k)$ in the word $u$ by $u(g_{1},\ldots ,g_{k})$.
\par
Let $u = v$ denote any semigroup identity. If $c(u) \cup c(v) =
\{x_{1},\ldots ,x_{k}\}$ let $w\in X^+ $ such that $c(w) = \{x_{1},\ldots ,x_{k}\}$ and
$| w| \ge 2$, and let $e \equiv w^{n-1}$.
We denote by $u^*=v^*$ the identity obtained from $u = v$ by
replacing each occurrence of $x_i$ by $ex_ie$.
\par
\vskip .3cm
\noindent{\bf LEMMA 1}. {\it Let ${\cal {V}}_n$ be the group variety determined by
the identity} (2) {\it and}
$$\{u_i= v_i, u_i,v_i\in X^*,i\in I\}.
\eqnum3$$
{\it Then ${\cal{SCRV}}_n$ is determined by
$$\{u^*_i= v^*_i, i\in I\}.
\eqnum{4}$$
and the identity} (1).
\par\vskip .3cm
\noindent{\bf PROOF}. Let ${\cal W}$ be the variety determined by (4) and
(1). Note that ${\cal W}\in L({\cal {SCR}}_n)$, since ${\cal W}$ satisfies
the identity (1). Consider any identity $u^*_i = v^*_i$ and let
$e_i \equiv w^{n-1}$, where $w\in X^+ $ such that
$| w| \ge 2 $ and $c(w) = c(u_i) \cup c(v_i) = \{x_1,\ldots,
x_k \}$. By Result 3 and its dual all elements $(e_ige_i)^{\tau }$,
where $g\in X^+$ and $c(g) \subseteq c(e_i)$, belong to the $\cal H$-class
of $e^{\tau }_i$. Since each $\cal H$-class in ${\cal{SCRV}}_n$-semigroup
is a group satisfying the identities of ${\cal {V}}_n$, this $\cal H$-class
satisfies the identity $u _i= v_i$ and we get
$$u^*_i\equiv u_i(e_ix_1e_i,\ldots , e_ix_ke_i)\;\tau\; v_i(e_ix_1e_i,\ldots
,e_ix_ke_i) \equiv v^*_i,$$
i.e., ${\cal{SCRV}}_n$ satisfies the identities (4). Hence ${\cal{SCRV}}_n \subseteq
 {\cal W}$. Conversely, consider any identity $u_i = v_i$ and, as
before, let $e_i\equiv w^{n-1}$, $w\in X^+ $ such that
$| w| \ge 2 $ and $c(w) = c(u_i) \cup c(v_i) = \{x_1,\ldots,
x_k \}$. Suppose $g_1, \ldots
, g_{k}\in X^+$ are
any words belonging to the same $\cal H$-class, when considered
as elements of $F{\cal W}$. Then the image of $\overline{e}_i\equiv e_i(g_1, \ldots
, g_{k})$ belongs to
the same $\cal H$-class, moreover it is the identity element of this
$\cal H$-class, i.e., $\overline{e}_ig_{j}\overline{e}_i \;\rho _{{\cal W}}\; g_{j}$
 for each~$j$. Hence
$$u_i(g_1,\ldots
, g_{k}) \;\rho _{\cal W}\; u_i(\overline{e}_ig_1\overline{e}_i,\ldots
,\overline{e}_ig_{k}\overline{e}_i) \equiv u^*_i(g_1,\ldots
,g_{k})
\;\rho _{\cal W}\; v^*_i(g_1,\ldots
,g_{k}) $$
$$
\equiv v_i(\overline{e}_ig_1\overline{e}_i,\ldots
,\overline{e}_ig_{k}\overline{e}_i) \;\rho _{\cal W}\; v_i(g_1,\ldots
,g_{k}).$$
So the identity $u_i= v_i$ holds in any $ \cal H$-class of $F{\cal W}$ and,
by the definition  of ${\cal{SCRV}}_n$, we get ${\cal W} \subseteq {\cal{SCRV}}_n$. Thus
${\cal{SCRV}}_n = {\cal W}$. \blacksquare
\vskip 0.2cm
\par
Before we are going to give a solution to the word problem for
$F{\cal{SCRV}}_n$ we will show  that it is not possible to apply the above mentioned  criterion
to the completely regular part of $T$ and $F$. It follows from
\par
\vskip .3cm
\noindent{\bf PROPOSITION 1}. {\it $T^2$ is not a free ${\cal CRV}_n$-semigroup. }
\par\vskip .3cm
\noindent{\bf PROOF}. First of all we remark that the ${\cal D}$-relation on
$T^2$ coincides with the restriction of the ${\cal D}$-relation
of $T$ on $T^2$. Indeed, let
$u,v\in T^2$ and $u{\cal D}v$ considered as elements
of $T$, i.e. $u=avb$ and $v=cud$ for some $a,b,c,d\in T$.
Then using identity (1) we have $u=(avb)^n=\overline av\overline b$
and $v=(cud)^n=\overline cu\overline d$ where $\overline a,\overline b,
\overline c,\overline d\in T^2$.
\par
Now suppose that $A=\{a_1,\ldots,a_p\}$ is the set of free generators
of $T^2$ considered as the free ${\cal {CRV}}_n$-semigroup.
Take the words $xy$ and $yx$. They belong to $T^2$ hence
$xy=a_k\ldots a_l$ and $yx=a_m\ldots a_s$, where $c(a_i)\subseteq \{x,y\}$ . By criterion of
equality of words in ${\cal {CRV}}_n$ $xy$ and $a_k\ldots a_l$  have the same left indicator in
$T$ by Result 3. Since $ a_k \neq x$  we have
$a_k=xyu$ for some $u \in T$ such that
$c(u)\subseteq \{x,y\}$. Similar $a_m=yxv$ for some $v \in T$ such that
$c(v)\subseteq \{x,y\}$. We see $a_k{\cal D}a_m$ considered as
elements of $T$ and as elements of $T^2$ too.
This implies $a_k\equiv a_m$ in $T^2$ but this is impossible
since $xyu$ and $yxv$ are not $\rho_{{\cal {CRV}}_n}$-related.
\blacksquare
\vskip 0.2cm
\par
Now we recall some definitions given in [5] or [7].
A {\it cell\/} of a word $u\in X^+$ is an ordered triple $(a,b,c)$ where
$a,b,c\in X^*$ such that $u \equiv abc$ and $b$ is a longest subword of
$u$ involving exactly $| c(u)| - 1$ letters. This means that if
$a$ is non-empty, then $a \equiv ex$ and if $c$ is non-empty, then
$c \equiv xf$ where $x\in X, e,f\in X^*$ and $x\not\in c(b)$. We call $a, b, c$ the
{\it left, middle\/} and {\it right\/} part of the cell $(a,b,c)$
respectively. We order all cells of a word $u$ putting
$(a_{0},b_{0},c_{0}), \ldots
 ,(a_{k},b_{k},c_{k})$ if
$| a_{0}| <| a_1| <\ldots
<| a_{k}| \; ($or $| c_{0}| >| c_1| >\ldots
>| c_{k}| )$.
Then the {\it characteristic sequence} of the word $u$ is the sequence
$(b_{0}, \ldots ,b_{k})$ of the middle parts of all cells ordered this way.
 We denote this characteristic sequence by $[u]$.
Note that then $b_{0}\equiv s(u)$ and $b_{k}\equiv d(u)$.
In addition we use the notation
$$
/w/_{u}=
\left\{\matrix{
[w]\hbox{ if }c(w)=c(u)\cr
w\hbox{ if }c(w)\subset c(u)\cr
}\right.
$$
\noindent and
$$
A_{u^{\rho_{\cal A}}} = \{w^{\rho_{\cal A}}: w\in X^+,c(w) \subset c(u), | c(w)| = |
c(u)| -1 \}
$$
for any variety $\cal A$.
\par
\vskip .3cm
\noindent{\bf DECOMPOSITION LEMMA}. {\it Let} $w\in X^+$ {\it and} $| c(w)| > 1.$
{\it Consider its characteristic sequence} $(w_{0},w_1,\ldots,w_{k})$ {\it and
let} $\{z_i\} = c(w)\backslash c(w_i)$ {\it for} $i = 0,\ldots,k.$
{\it There exists} $d_i$ {\it such that} $d_i\equiv d(w_{i-1}) \equiv s(w_i)$
{\it for} $ i = 1,\ldots,k.$ {\it Then for any words} $f_i,g_i\in X^*,
i = 1,\ldots,k$ {\it such that} $c(f_i), c(g_i) \subseteq c(w)$
$$
w\;\rho \; w_{0}z_{0}g_1(f_1z_1d_1z_{0}g_1)^{n-2}f_1z_1w_1z_1g_2
(f_2z_2d_2z_1g_2)^{n-2}f_2z_2w_2\ldots f_{k}z_{k}w_{k}
$$
\noindent {\it holds.}
\par\vskip .3cm
\noindent{\bf PROOF}. If $w\in X^+,| c(w)| >1$ and $(a_{0},w_{0},c_{0}),\ldots
,(a_{k},w_{k},c_{k})$ is the sequence of all cells of the word $w$
ordered by $| a_{0}| <\ldots
<| a_{k}| $ and $\{z_i\} = c(w)\backslash c(w_i)$ for
$i = 0,\ldots ,k$, then there are $d_i,e_i,h_i\in X^*$ such that $w_{i-1}\equiv
e_iz_id_i, w_i\equiv d_iz_{i-1}h_i$ and $w \equiv
a_{i-1}e_iz_id_iz_{i-1}h_ic_i$. Note
that $d_i \equiv d(w_{i-1}) \equiv s(w_i)$. If $f_i,g_i\in X^*, i = 1,\ldots
,k$ are
any words such that $c(f_i),c(g_i) \subseteq c(w)$, then by Result 4
(taking $g_if_i$ as $v, z_id_iz_{i-1}$ as $u$ and since $c(z_id_iz_{i-1}) =
c(w)$ and $z_id_iz_{i-1}$ is a group element when considered
as an element of $F{\cal {SCR}}_n)$ we get
$$
w_{i-1}z_{i-1}h_i \equiv e_iz_id_iz_{i-1}h_i \;\rho
$$
$$
e_iz_id_iz_{i-1}(g_if_iz_id_iz_{i-1})^{n-1}h_i \equiv
w_{i-1}z_{i-1}g_i(f_iz_id_iz_{i-1}g_i)^{n-2} f_iz_iw_i.
$$
\noindent Then taking into account that $c_{i-1}\equiv z_{i-1}h_ic_i$, we can
write the sequence
$$
w \equiv w_{0}c_{0}\equiv w_{0}z_{0}h_1c_1 \;\rho
\;w_{0}z_{0}g_1(f_1z_1d_1z_{0}g_1)^{n-2}f_1z_1w_1c_1 \;\rho \ldots
$$
$$
\rho\; w_{0}z_{0}g_1(f_1z_1d_1z_{0}g_1)^{n-2}f_1z_1w_1z_1g_2
(f_2z_2d_2z_1g_2)^{n-2}f_2z_2w_2z_2 \ldots
 f_{k}z_{k}w_{k}.\hbox{ \blacksquare
}$$
\par
\vskip .2cm
Now we are ready to give a solution to the word
problem for $F{\cal {SCR}}_n$. If $c(u) = c(v) = \{x\}$, i.e. $u \equiv x^{m}, v \equiv x^{k}$
for some
$m,k$, then $u \;\rho\; v$ if and only if either $m = k = 1$ or $m,k>1$
and $n-1$ divides $m-k$.
\par
The following theorem describes the case $| c(u)| > 1$.
\par
\vskip .3cm
\noindent{\bf THEOREM 1}. {\it Let} $u,v \in X^+, | c(u)| > 1$ . {\it Then}
 $u\;\rho\; v$ {\it if and only if }
\par
\noindent (1) $c(u) = c(v)$,
\par
\noindent (2) $s(u)\; \rho \;s(v)$,
\par
\noindent (3) $d(u) \;\rho \; d(v)$,
\par
\noindent {\it and}
\par
\noindent (4) {\it if} $(u_{0}, \ldots ,u_{k})$
{\it and} $(v_{0}, \ldots ,v_{l})$ {\it are the
characteristic sequences of} $u$ {\it and} $v$ {\it respectively and if}
$\varphi : A_{u^\rho} \rightarrow X$ {\it is any injective mapping then}
\par
$$
\varphi (u_1^\rho) \ldots
 \varphi (u_{k-1}^\rho) \;\rho _{{\cal G}_n} \varphi (v_1^\rho) \ldots
 \varphi (v_{l-1}^\rho)\hbox{ {\it in} }X^*.
$$
\par\vskip .3cm
\noindent{\bf PROOF}. {\sl Necessity.} If $u \;\rho \;v$ then obviously $c(u)
= c(v)$ and the conditions (2) and (3) are true by
Result 1 and its dual.
\par
We only need to show the necessity of (4) in the case when $u \equiv pqh$
and $v \equiv pq^nh$, where $| q| \ge 2.$
Two cases are possible:
\par
Case 1.$\;   c(q) \subset c(u)$. Let $(u_{0}, \ldots ,u_{k})$ be the
characteristic sequence of $u \equiv pqh$. For each $0\le i \le k$, let
$(a_i, u_i, b_i)$ be the cell which corresponds to the middle part $u_i$.
For each such cell we shall find a corresponding cell $(c_i, v_i, d_i)$ of the word
$v \equiv pq^nh$ such that $u_i \;\rho\; v_i$. In view of the fact that $c(q)\subseteq c(u)$,
$u_i$ cannot be a proper subword of $q$. Therefore only the following cases can occur.
\par
1. $ | a_iu_i | < | pq |$. In this case $pq \equiv a_iu_ig_i$ for some
$g_i \in  X^+$ and  $(c_i, v_i, d_i)$ is a corresponding cell for $v$ with
$c_i \equiv a_i$, $ v_i \equiv u_i$ and $d_i \equiv g_iq^{(n-1)}b_i$.
\par
2. $ |a_iu_i | \ge | pq |$, $| a_i | \le | p |$.
In this case $u_i \equiv g_iqh_i$ for some $g_i , h_i\in  X^*$ and
$(c_i, v_i, d_i)$ is a corresponding cell for $v$ with
$c_i \equiv a_i$, $ v_i \equiv g_iq^nh_i$ and $d_i \equiv b_i$.
We remark that in this case $ u_i \equiv g_iqh_i \;\rho\; g_iq^nh_i \equiv v_i $ .
\par
3. $ |a_iu_i | \ge | pq |$, $| a_i | > | p |$.
In this case $qh \equiv g_iu_ib_i$ for some $g_i \in  X^+$ and
$(c_i, v_i, d_i)$ is a corresponding cell for $v$ with
$c_i \equiv a_iq^{(n-1)}g_i$, $ v_i \equiv u_i$ and $d_i \equiv b_i$.
\par
It is now easy to see that the sequence $(v_0, \ldots ,v_k) $ constructed in this way
is the characteristic sequence of $v$ and that in fact
$$
\varphi (u_1^\rho) \ldots
 \varphi (u_{k-1}^\rho) \equiv \varphi (v_1^\rho) \ldots
 \varphi (v_{k-1}^\rho).
$$
\par
Case 2. If $c(q)  = c(u)$ then the situation where $q$
belongs to a middle part for some cell is impossible. It
is however possible that $u$ will have cells whose middle
part will be a subword of $q$. Let $[q] = (q_{0},\ldots
,q_{m})$ be the
characteristic sequence of $q$. We begin to construct the
characteristic sequence of $u$. First of all, it contains
either $pq_{0}$ if $c(pq_{0}) \subset c(u)$ or $[pq_{0}]$ otherwise, then $q_1,\ldots
,q_{m-1}$, and finally either $q_{m}h$ or $[q_{m}h]$,
i.e.
$$
[u] = (/pq_{0}/_{u},q_1,\ldots
,q_{m-1},/q_{m}h/_{u})
$$
\noindent By analogy we get
$$
[v] = (/pq_{0}/_{u},(q_1,\ldots
,q_{m-1},/q_{m}q_{0}/_{u})^{n-1},
q_1,\ldots
,q_{m-1},/q_{m}h/_{u})
$$
\noindent Thus we see that in $F{\cal G}_n$ the images of these
sequences are equal.
\par
{\sl Sufficiency}. Let $u,v\in X^+$ with
$| c(u)| = | c(v)|  > 1$ such that
the conditions (1),(2),(3) and (4) of the statement of
the Theorem are satisfied. Let $c(u) = c(v) =
\{x_1,\ldots ,x_m\}$. For $w \in X^+$ with
$c(w) \subseteq \{x_1,\ldots ,x_m\}$,
$| c(w)| = m-1$ we define:
$$
\alpha (w) \equiv x_1\ldots
x_mz(w)wz(w)x_1\ldots
x_m(x_1\ldots
x_m\delta(w)d(w)z(w)x_1\ldots
x_m)^{n-2},
$$
$$
\alpha _{L}(w) \equiv wz(w)x_1\ldots
x_m(x_1\ldots
x_m\delta(w)d(w)z(w)x_1\ldots
x_m)^{n-2},
$$
$$
\alpha _{R}(w) \equiv x_1\ldots
x_mz(w)w,
$$
where $\{z(w)\} = \{x_1,\ldots
,x_m\}\backslash c(w)$. Using the Decomposition Lemma and taking
$f_i\equiv g_i\equiv x_1\ldots
x_m$
for all $i$ we get
$$
u \;\rho \;\alpha _{L}(u_{0})\alpha (u_1)\ldots
\alpha (u_{k-1})\alpha _{R}(u_{k})
$$
\noindent and
$$
v \;\rho \;\alpha _{L}(v_{0})\alpha (v_1)\ldots
\alpha (v_{l-1})\alpha _{R}(v_{l}).
$$
\par
Consider the condition (4). For any
$a, b\in \{u_1,\ldots ,u_{k-1},v_1,\ldots ,v_{l-1}\}$
and any map $\varphi : A_{u^\rho} \rightarrow X $, $ \varphi (a^\rho)\equiv
\varphi (b^\rho)$
implies $\alpha (a)\;\rho\; \alpha (b)$ since then $a \;\rho \;b$ and so $z(a) \equiv z(b)$ ,
$d(a) \;\rho \;d(b)$  and $\delta(a) \equiv \delta(b)$ by the necessity part of the proof.
 By definition $ \alpha (u_1),\ldots
,\alpha (u_{k-1}),
\alpha (v_1),\ldots
,\alpha (v_{l-1})$, when considered as elements of $F$, belong to the same
$\cal H$-class, hence if
$$\varphi (u_1^\rho)\ldots
\varphi (u_{k-1}^\rho) \;\rho _{{\cal G}_n} \;\varphi (v_1^\rho)\ldots
\varphi (v_{l-1}^\rho),$$ then
$$\alpha (u_1)\ldots
\alpha (u_{k-1}) \;\rho \;\alpha (v_1)\ldots
\alpha (v_{l-1}).$$
\par
From the conditions (2),(3) we see also that
$$
\alpha _{L}(u_{0}) \equiv \alpha _{L}(s(u)) \;\rho \;\alpha _{L}(s(v)) \equiv \alpha _{L}(v_{0})
$$
\noindent and
$$
\alpha _{R}(u_{k}) \equiv \alpha _{R}(d(u)) \;\rho\; \alpha _{R}(d(v)) \equiv \alpha _{R}(v_{l})
$$
\noindent too. Thus we have $u \;\rho\; v$. \blacksquare
\par
\vskip .2cm
As a solution to the word problem for the free
${\cal{SCRV}}_n$-semigroup on $X$ we will give an algorithm reducing this
problem to the word problem for the free group on $X$ of
the group variety ${\cal {V}}_n$. In order to prove this criterion we need the next
\par
\vskip .3cm
\noindent{\bf LEMMA 2}. {\it Let} $u\in X^+,| c(u)| > 1, [u] = (u_{0},\ldots
,u_{k}),
/d(u)s(u)/_{u}= (\overline u_{0},\ldots, \overline u_{m})$({\it if} $c(d(u)s(u))\subset c(u) ,$
{\it then} $m=0,  \overline u_0=d(u)s(u))$ {\it and}
$$\psi : A_{u^\rho} \rightarrow X$$
{\it an injective
mapping. Let } $\overline{u}= \psi (\overline u_{0}^{\rho})
\ldots\psi (\overline u_{m}^{\rho })${\it . Then the mapping}
$\overline{\psi}$ {\it from the} $\cal H$-{\it class of} $u^{\rho }$ {\it in}
$F{\cal {SCR}}_n$ {\it into} $F{\cal G}_n${\it , given by }
$$
\overline{\psi}(u^{\rho }) = \psi (u_1^{\rho })
\ldots\psi (u_{k-1}^{\rho })\overline{u}
$$
{\it is an injective homomorphism.}
\par\vskip .3cm
\noindent{\bf PROOF.} Let $u^{\rho}{\cal H} v^{\rho}$
and $[v] = (v_{0},\ldots,v_{l})$.
From Result 3 and its dual we have that $d(u)s(u) \;\rho\; d(u)s(v) \;\rho\; d(v)s(v)$, and
thus in particular, $c(d(u)s(u)) = c(d(u)s(v)) = c(d(v)s(v))$.
If $c(d(u)s(u))\subset c(u)$, then $\overline u \equiv \psi(d(u)s(u)^{\rho})\equiv
\psi(d(v)s(v)^{\rho})
\equiv \overline v.$ If  $c(d(u)s(u)) = c(u)$ then put
$$
[d(u)s(v)] = (\overline w_{0},\ldots, \overline w_{s}),
$$
$$
[d(v)s(v)] = (\overline v_{0},\ldots, \overline v_{t}).
$$
A proof of the fact that $d(u)s(u)$ and $d(u)s(v)$ are $\rho$-related makes use only
of transitions as considered in Case 1 of  the necessity part of the proof of Theorem 1.
As in the proof of Theorem 1 we may conclude here that $m=s$ and
$\psi (\overline u_i) = \psi (\overline w_i)$ for every $0\le i\le m$. A similar argument yields
$ s=t$ and $\psi (\overline w_i) = \psi (\overline v_i)$ for every $0\le i\le s$.
In conclusion, $m=t$, $\psi (\overline u_i) = \psi (\overline w_i) = \psi(\overline v_i)$ for every
$0\le i\le m$, and so
$\overline u = \psi (w_{0}^{\rho}) \ldots\psi (w_{m}^{\rho }) =\overline v$.
\par
Again let $v^{\rho}{\cal H} u^{\rho }$ . Since
$$
[uv] = (u_{0},\ldots
,u_{k-1}, /d(u) s(v)/_{u},v_1,\ldots,v_{l})
$$
we have
$$
\overline{\psi}(uv^{\rho }) = \psi (u_1^{\rho })\ldots
\psi (u_{k-1}^{\rho })\overline{u}\psi (v_1^{\rho})
\ldots\psi (v_{l-1}^{\rho })\overline{uv}
$$
$$
= \psi (u_1^{\rho })\ldots \psi (u_{k-1}^{\rho })
\overline{u}\psi (v_1^{\rho})\ldots\psi (v_{l-1}^{\rho })
\overline{v}=\overline{\psi}(u^{\rho})
\overline{\psi}(v^{\rho}),$$
where $ \overline {uv}= \overline v$ follows from $u^{\rho}{\cal H} uv^{\rho} $.
That $\overline \psi $ is injective follows immediately from Theorem 1.
\blacksquare
\vskip .2cm
\par
The following theorem gives a solution to the word
problem for $F{\cal{SCRV}}_n.$
\par
\vskip .3cm
\noindent{\bf THEOREM 2}. {\it Let} $u,v \in X^+, | c(u)| > 1.$ {\it Then}
 $u\;\tau \; v$ {\it if and only if}
\par
\noindent (1) $c(u) = c(v)$,
\par
\noindent (2) $s(u) \;\tau \;s(v)$,
\par
\noindent (3) $d(u) \;\tau \;d(v)$,
\par
\noindent {\it and}
\par
\noindent (4) {\it if} $(u_{0}, \ldots
 ,u_{k})$ {\it and} $(v_{0}, \ldots
 ,v_{l})$ {\it are the
characteristic sequences of} $u$ {\it and} $v$ {\it respectively and if}
$\varphi : A_{u^\tau} \rightarrow X$ {\it is any injective mapping then}
\par
$$
\varphi (u_1^\tau) \ldots
 \varphi (u_{k-1}^\tau) \;\rho _{{\cal {V}}_n} \;\varphi (v_1^\tau) \ldots
 \varphi (v_{l-1}^\tau)\hbox{ {\it in} }X^+.
$$
\par\vskip .3cm
\noindent{\bf PROOF}. {\sl Necessity.} Let $d,f \in X^+$ such that $d \;\rho_{{\cal V}_n}\; f$, that is,
$d=f$ is a semigroup identity which holds true in the group variety $\cal V_n$.
Let $c(d) \cup c(f) = \{x_1,\ldots ,x_{k}\}$. We can write the identity $d=f$
also as $d(x_1, \ldots , x_{k}) = f(x_1, \ldots , x_{k})$. Let $w \in X^+$ such that
$c(w) = \{x_1,\ldots ,x_k\}$ and $| w | \ge 2$, and put  $e = w^{n-1}$. A
substitution of $x_i$ by $ex_ie$ yields the identity
$$
d(ex_1e, \ldots , ex_ke) = f(ex_1e,  \ldots ,
ex_ke).
$$
The set of all the pairs $(d(ex_1e, \ldots , ex_ke) , f(ex_1e,  \ldots , ex_ke))$
thus obtained will be denoted by $\mu $. From Lemma 1 (and its proof) it follows that
$\tau$ is the fully invariant congruence relation generated by $\rho \cup\mu$.
\par
It is easy to see that $\mu$ is a fully invariant equivalence relation, hence the congruence
 relation generated by $\mu$ is fully invariant. Therefore $\tau$ is the transitive closure of the
union of  $\rho$ with  the congruence generated by $\mu$. It therefore suffices to show that
the conditions (1), (2), (3) and (4) are satisfied  if either  $u\;\rho\; v$ or $u\equiv
ad(ex_1e, \ldots , ex_ke)b,\; v\equiv af(ex_1e,  \ldots , ex_ke)b$ with $d\;\rho_{{\cal V}_n}\; f$,
$e$ as given above, and $a,b \in X^*$.
\par
That (1) holds is obvious from Theorem 1. We prove that (2) is satisfied: the proof for (3)
follows dually. If $u\;\rho \; v$ then then $s(u)\;\rho\; s(v)$ by Theorem 1, and thus
$s(u)\;\tau\; s(v)$ since $\rho \subseteq \tau$. Otherwise $u\equiv
ad(ex_1e, \ldots , ex_ke)b, \; v\equiv af(ex_1e,  \ldots , ex_ke)b$
and either one of the following cases occur:
\par
1.$ \; s(u) \equiv s(a)\equiv s(v)$
\par
2.$ \; s(u) \equiv ap\equiv s(v)$ for some initial segment $p$ of $e$
\par
\noindent or
\par
3.$\; s(u) \equiv ad(ex_1e, \ldots , ex_ke)p,\; s(v)\equiv af(ex_1e,  \ldots , ex_ke)p$ for some
initial segment $p$ of $b$.
\par
\noindent
In the last case $s(u)$ and $s(v)$ are related in the congruence relation generated by $\mu$. Thus, in all cases, $s(u)\;\tau\; s(v)$.
\par
We  next show  that (4) holds. Clearly if $u\;\rho \; v$ then (4) holds by Theorem~1.
Otherwise, we have $u\equiv
ad(ex_1e, \ldots , ex_ke)b, \; v\equiv af(ex_1e,  \ldots , ex_ke)b$.
We shall put $d^* = d(ex_1e, \ldots , ex_ke)$ and $f^* = f(ex_1e,  \ldots , ex_ke)$.
We first consider the case where $c(a)\cup c(b) \subseteq \{x_1,\ldots ,x_k\}.$
Let
$$
[d^*] = (d_0,\ldots , d_s), [f^*] = (f_0,\ldots , f_t).
$$
We have $ as(d^*) \equiv as(e) \equiv as(f^*)$ and $d(d^*)b \equiv d(e)b \equiv d(f^*)b$, thus
$$
[u] = (a_0,\ldots , a_p, d_1,\ldots , d_{s-1}, b_0,\ldots b_q),
$$
$$
[v] = (a_0,\ldots , a_p, f_1,\ldots , f_{t-1}, b_0,\ldots b_q),
$$
where
$$
/as(d^*)/_u = /as(f^*)/_v = (a_0,\ldots , a_p),$$
$$
/d(d^*)b/_u = /d(f^*)b/_v = (b_0,\ldots , b_q).$$
\par
Let  $\psi : A_{u^{\rho}} \rightarrow X$ be given by
$\psi (c^\rho) = \varphi (c^\tau)$, and let $\overline \psi$ be the injective homomorpism
from the $\cal H$-class of $F$ containing $e$ into $F\cal G_n$. Note that this
$\cal H$-class contains $d^*, f^*$ and the $ex_ie, 1\le i\le k$. As in Lemma 2
we have
$$
\psi (d_1^\rho)\ldots \psi (d_{s-1}^\rho)\overline {d^*} = \overline \psi (d^*{}^\rho) =
d( \overline \psi (ex_1e^\rho),\ldots ,\overline \psi (ex_ke^\rho))
$$
$$
\rho_{{\cal V}_n}\; f(\overline \psi (ex_1e^\rho),\ldots ,\overline \psi (ex_ke^\rho)) =
\overline \psi (f^*{}^\rho) = \psi (f_1^\rho)\ldots \psi (f_{t-1}^\rho)\overline {f^*}
$$
with $\overline {d^*}  = \overline {f^*}$ in $F\cal G_n$. Therefore
$$
\varphi (d_1^\tau)\ldots \varphi (d_{s-1}^\tau) \;
\rho _{{\cal V}_n}\; \varphi (f_1^\tau)\ldots
\varphi (f_{t-1}^\tau)
$$
and thus also
$$
\varphi (a_1^\tau)\ldots
\varphi (a_{p}^\tau)\varphi (d_1^\tau)\ldots
\varphi (d_{s-1}^\tau)\varphi (b_0^\tau)\ldots \varphi (b_{p-1}^\tau)
$$
$$
\rho _{{\cal V}_n} \; \varphi (a_1^\tau)\ldots
\varphi (a_{p}^\tau)\varphi (f_1^\tau)\ldots
\varphi (f_{t-1}^\tau) \varphi (b_0^\tau)\ldots \varphi (b_{p-1}^\tau) .
$$
\par
If $c(a)\cup c(b)\not\subseteq  \{x_1,\ldots , x_{k}\}$,
then the middle parts of the cells for $u$ and $v$ may contain $d^*$ and $f^*$
respectively as subwords. In fact, $a'd^*b'$ is such an entry in $[u]$ if and only if
the corresponding $a'f^*b'$ is an entry in $[v]$. Since both $d^*$ and $f^*$ begin
and end with $e$, and $c(e) = c(d^*) = c(f^*)$, the other entries in $[u]$ and $[v]$
also correspond and such corresponding entries are identical in $X^+$.
In conclusion, the characteristic sequences for $u$ and $v$ are of equal length and corresponding
entries are $\tau$-related. Therefore in this case the condition (4) is trivially satisfied.
\par
{\sl Sufficiency}. An argument similar to the one used in the proof of  Theorem~1, replacing
$\rho$ with $\tau$ and $\rho_{{\cal G}_n}$ with $\rho_{{\cal {V}}_n}$ respectively, takes care of this
part.
\blacksquare
\par
\vskip .3cm
\centerline{{\bf 3. The join of ${\cal{CR}}_n$ and ${\cal{SI}}$}
}
\vskip .3cm
In order to answer the question whether ${\cal {SCR}}_n$
coincides with ${\cal{CR}}_n\vee {\cal{SI}}$ we shall first find a basis of
identities for this join. This will be done in the next
\par
\vskip .3cm
\noindent{\bf LEMMA 3}. ${\cal{CR}}_n\vee {\cal{SI}} = var \{(xy)^n = xy,\ xy^nx = xyx \}.$
\par\vskip .3cm
\noindent{\bf PROOF}. Clearly ${\cal{CR}}_n$ and ${\cal{SI}}$ satisfy the identity (1)
and ${\cal{CR}}_n$ satisfies the identity
$$xy^nx = xyx.
\eqnum{5}$$
The fact that ${\cal{SI}}$ satisfies identity (5) follows
from the solution to the word problem for ${\cal{SI}}$ given by
Gerhard in [2].
\par
Conversely, suppose $u = v$ is any identity which
holds in ${\cal{SI}}\vee {\cal{CR}}_n$. Consider the semigroup $S$ which is the
free object on a countably infinite number of variables of the variety given by the identities (1)
and (5). It suffices to prove that $u = v$ holds in $S$. We
certainly have that $c(u) = c(v)$. We use induction on the
number of letters in $c(u)$.
\par
First consider the case where $| c(u)| = 1.$ Any
identity $u = v$ which holds both in ${\cal{CR}}_n$ and ${\cal{SI}}$ is a
consequence of $x^2= x^{n+1}$, but clearly this identity holds
in $S$ too.
\par
Now let $u$ and $v$ be such that $c(u) = c(v)$ contains
$k > 1$ letters and assume that  an identity involving
less than $k$ letters which holds in ${\cal{CR}}_n\vee {\cal{SI}}$ holds in $S$
too. Consider $s(u)$ and $s(v), d(u)$ and $d(v)$. By the known
criterion of equality of words in ${\cal{CR}}_n$ (see [7] or [5])
and in ${\cal{SI}}, s(u) = s(v)$ in ${\cal{CR}}_n$ and in ${\cal{SI}}$,
and since $s(u)$ depends on $k-1$ letters $s(u) = s(v)$ in $S$.
Dually, $d(u) = d(v)$ in $S$. Moreover $\sigma (u) \equiv \sigma (v)$
and $\delta(u) \equiv \delta(v)$, hence
$s(u)\sigma (u) = s(v)\sigma (v)$ and $\delta (u)d(u) =\delta (v)d(v)$ in
$S$. Then
we use the notation $s(u)\sigma (u) \equiv a$ and $\delta (u)d(u) \equiv b$ and
we consider a sequence of applications of the identity
$x^n= x$ transforming $u$ into $v$:
$$
u \equiv u_{0}\rightarrow u_1\rightarrow \ldots
\rightarrow u_{k}\equiv v,
$$
\noindent where
$u_i \equiv pw_iq\hbox{ and }u_{i+1} \equiv pw^n_iq\hbox{ or vice versa.
}$
 Multiplying by $a^{n-1}$ from the left and by $b^{n-1}$ from the
right and taking into consideration that $a^n= a$ and $b^n= b$
in $S$ (since they depend on more than one letter) we obtain
$u = a^{n-1}u_0b^{n-1},\; v = a^{n-1}u_kb^{n-1}$ in $F{\cal SCR}_n$ and
$$
a^{n-1}u_{0}b^{n-1} \rightarrow a^{n-1}u_1b^{n-1} \rightarrow \ldots
 \rightarrow a^{n-1}u_{k}b^{n-1}
$$
\noindent where
$a^{n-1}u_ib^{n-1} \equiv p_1w_iq_1$ and $a^{n-1}u_{i+1}b^{n-1} \equiv p_1w^n_iq_1$ and
where
$c(p_1) = c(q_1)$ since $c(a)=c(b)=c(u)$. Then using Result 4 taking firstly
$p_1$ and
$q_1$ as $u$ and $v$ and then $q_1$ and $p_1$ as $u$ and $v$, we obtain:
$$
p_1w_iq_1 = p_1(q_1p_1)^{n-1}w_i(q_1p_1)^{n-1}q_1
$$
$$
p_1w^n_iq_1 = p_1(q_1p_1)^{n-1}w^n_i(q_1p_1)^{n-1}q_1
$$
\noindent in $F{\cal {SCR}}_n$, i.e., we have a sequence transforming $u$ into $v$
in which every step may be considered as an application
of the identity (5). Thus we see that $u = v$ holds in $S$.\blacksquare
\par
\vskip .2cm
Now the above proved Lemma and Theorem 1 enable us
to obtain
\par
\vskip .3cm
\noindent{\bf THEOREM 3}. ${\cal {SCR}}_n \neq {\cal{CR}}_n \vee {\cal{SI}}.$
\par\vskip .3cm
\noindent{\bf PROOF}. We prove that the identity (5) does not hold
in ${\cal {SCR}}_n$.
The characteristic sequences of $xyx$ and $xy^nx$ are :
$$
[xyx] = (x,y,x),\; [xy^nx] = (x,y^n,x).
$$
Since $y \neq y^n$ in $F{\cal {SCR}}_n$ we have
$xyx \neq xy^nx$ in $F{\cal {SCR}}_n$ by Theorem 1.
\blacksquare
\par
\vskip .3cm
\centerline{{\bf 4. The semigroups with orthogroup square}
}
\vskip .3cm
In this section we shall consider only the variety
${\cal{SO}}_n$. The variety under consideration,
in contrast with the situation for the variety of all
${\cal {SCR}}_n$-semigroups, is the join of ${\cal O}_n$ and ${\cal{SI}}$. In order to
prove this fact we need
\par
\vskip .3cm
\noindent{\bf LEMMA 4}. {\it The identity}
\par
$$
 x^{2n-2}yx^{2n-2}= x^{2n-2}y^nx^{2n-2}
$$
\noindent {\it holds in} $F{\cal{SO}}_n$.
\par\vskip .3cm
\noindent{\bf PROOF}. Firstly we show that $x^{2n-2}y^{k}x^{2n-2} = (x^{2n-2}yx^{2n-2})^{k}$
holds in
$F{\cal{SO}}_n$ for any $k\le n$. We use induction on $k$. The case $k = 1$ is
obvious. Assume that $k > 1$ and take two idempotents
$(x^{2n-2}y^{k-1})^{n-1}$ and $(yx^{2n-2})^{n-1}$. Since $F{\cal{SO}}_n^2$ is orthodox,
their
product $(x^{2n-2}y^{k-1})^{n-1} (yx^{2n-2})^{n-1}$ must be an idempotent,
belonging to the $\cal H$-class of $(x^{2n-2}y^{k-1})^{n-1}(yx^{2n-2})^{n-1}$
considered as element of $F{\cal{SO}}_n$.
From Theorem 1 we have that $x^{2n-2}yx^{2n-2}$ and
$(x^{2n-2}y^{k-1})^{n-1}(yx^{2n-2})^{n-1}$ are $\cal H$-related when
 considered as elements of
$F{\cal {SCR}}_n$ and hence as elements of $F{\cal{SO}}_n$ too. We obtain that
in $F\cal {SO}_n$
$$
(x^{2n-2}y^{k-1})^{n-1} (yx^{2n-2})^{n-1}= (x^{2n-2}yx^{2n-2})^{n-1},
$$
$$
(x^{2n-2}y^{k-1}x^{2n-2})^{n-2} x^{2n-2}y^{k}x^{2n-2} (x^{2n-2}yx^{2n-2})^{n-2}=
(x^{2n-2}yx^{2n-2})^{n-1}.
$$
Using the induction assumption and the identity (1) we
can write
$$
(x^{2n-2}y^{k-1}x^{2n-2})^{n-2} = (x^{2n-2}yx^{2n-2})^{(n-2)(k-1)}
$$
$$ = (x^{2n-2}yx^{2n-2})^{nk-2k-n+2} =
 (x^{2n-2}yx^{2n-2})^{n(k-2)+n-2k+2}
$$
$$
= (x^{2n-2}yx^{2n-2})^{k-2+n-2k+2} = (x^{2n-2}yx^{2n-2})^{n-k}.
$$
\noindent Thus we get
$$
(x^{2n-2}yx^{2n-2})^{n-k} x^{2n-2}y^{k}x^{2n-2} (x^{2n-2}yx^{2n-2})^{n-2} = (x^{2n-
2}yx^{2n-2})^{n-1}
$$
\noindent Multiplying this equality by $(x^{2n-2}yx^{2n-2})^{k-1}$ from the left and
by $x^{2n-2}yx^{2n-2}$ from the right we get
$$
(x^{2n-2}yx^{2n-2})^{n-1}x^{2n-2}y^{k}x^{2n-2}(x^{2n-2}yx^{2n-2})^{n-1} = (x^{2n-
2}yx^{2n-2})^{k}.
$$
\noindent Taking into account that $(x^{2n-2}yx^{2n-2})^{n-1}$ is the identity of
this $\cal H$-class we obtain
$$
x^{2n-2}y^{k}x^{2n-2} = (x^{2n-2}yx^{2n-2})^{k}.
$$
\noindent Finally we take $k = n$ and obtain
$$x^{2n-2}y^nx^{2n-2} = (x^{2n-2}yx^{2n-2})^n =
x^{2n-
2}yx^{2n-2}.\eqno{ \blacksquare
}$$
\par
\vskip .2cm
Before we state the main result of this section we
make the following remark. If we consider ${\cal O}_n$ instead of
${\cal{CR}}_n$ in Lemma 3 then we can prove
\par
\vskip .3cm
\noindent{\bf LEMMA 3}$^*.$ $ {\cal{SI}}\vee {\cal O}_n = var \{(xy)^n = xy,xy^nx =
xyx,
((xy)^{n-1}(zt)^{n-1})^2 = (xy)^{n-1}(zt)^{n-1}\}.
$
\blacksquare
\par
\vskip .2cm
\noindent This fact holds because the varieties ${\cal{SI}}$ and ${\cal O}_n$
are such that the squares of their semigroups are
orthodox.
\par
\medskip
\noindent{\bf THEOREM 4}. ${\cal{SI}}\vee {\cal O}_n = {\cal{SO}}_n$.
\par\vskip .3cm
\noindent{\bf PROOF}. By the previous Lemma we have that ${\cal{SI}}\vee {\cal O}_n \subseteq
{\cal{SO}}_n$
because the identity (1) obviously holds in ${\cal{SI}}\vee {\cal O}_n$ and the
square of an orthodox semigroup is orthodox. It remains to
verify that ${\cal{SO}}_n \subseteq {\cal{SI}}\vee {\cal O}_n$,
i.e. the identity $xy^nx = xyx$
holds in ${\cal{SO}}_n$. Using (1) we get $xy^nx = (xy)^{n-1}xy^nx(yx)^{n-1}$.
We put $a \equiv (xy)^{n-1}x$ and $b \equiv x(yx)^{n-1}$. Then we apply
Lemma 4 and Result 4, taking firstly $a, b$ for $u, v$ and
then for $v, u$
$$
xy^nx = ay^nb = a(ba)^{n-1} y^n(ba)^{n-1}b = a(ba)^{n-1} y(ba)^{n-1}b =
ayb = xyx .\eqno{ \blacksquare}
$$
\par
\noindent{\bf COROLLARY 1}. {\it The lattice} $L({\cal{SO}}_n)$ {\it is a subdirect
product of the lattice} $L({\cal O}_n)$ {\it and of the interval}
$[{\cal{SI}} \cap {\cal O}_n, {\cal{SI}}]$.
\par\vskip .3cm
\noindent{\bf PROOF}. By Theorem III.2.4.[4] $L({\cal{SO}}_n)$ is a
subdirect product of the $L({\cal O}_n)$ and $[{\cal O}_n, {\cal{SO}}_n]$. But
$L({\cal{SO}}_n)$ is
modular as sublattice of $L({\cal {SCR}}_n)$ (by Corollary 3.3 [11]),
whence by Theorem IV.1.2.[4] $[{\cal O}_n, {\cal O}_n\vee {\cal{SI}}] \cong
[{\cal{SI}} \cap {\cal O}_n, {\cal{SI}}]$.
\phantom{text}\blacksquare
\par
This result is important since in [8] and [9] the
lattice of all varieties of orthogroups was investigated
and a subdirect decomposition of it was found. The same
was done by Gerhard in [3] for the lattice of all
varieties of semigroups with idempotent square, moreover
it was shown to be distributive. Thus we get the full
description of the lattice $L({\cal{SO}}_n)$.

\end{document}